\documentclass[11pt]{amsart}
\usepackage{amssymb,amsmath,amsfonts,amscd,euscript}

\newcommand{\nc}{\newcommand}

\nc{\slt}{\mathfrak{sl}_2}
\nc{\suth}{\widehat{\mathfrak{su}}(2)}
\nc{\gl}{\mathfrak{gl}}
\nc{\GL}{\mathfrak{GL}}
\nc{\g}{\mathfrak{g}}
\nc{\gh}{\widehat{\mathfrak{g}}}
\nc{\h}{\mathfrak{h}}
\nc{\hh}{\widehat{\mathfrak{h}}}
\nc{\la}{\lambda}
\nc{\slth}{\widehat{\slt}}
\nc{\C}{\mathbb C }
\nc{\Z}{\mathbb Z }
\nc{\N}{\mathbb N }
\nc{\al}{\alpha}
\nc{\be}{\beta}
\nc{\ve}{\varepsilon}
\nc{\ch}{{\mathop {\rm ch}}}
\nc{\Id}{{\mathop {\rm Id}}}
\nc{\Tr}{{\mathop {\rm Tr}\,}}
\nc{\U}{{\mathop {\rm U}}}
\nc{\bra}{\langle}
\nc{\ket}{\rangle}
\nc{\ld}{\ldots}
\nc{\cd}{\cdots}
\nc{\hk}{\hookrightarrow}
\nc{\n}{\mathfrak{n}}
\nc{\un}{\mathfrak{u}}
\nc{\bo}{\mathfrak{b}}
\nc{\T}{\otimes}
\nc{\wt}{\widetilde}

\nc{\qb}[2]{\genfrac{(}{)}{0pt}{}{#1}{#2}_q}
\nc{\fac}[1]{(#1)_q!}

\newtheorem{theo}{Theorem}[section]
\newtheorem*{theo*}{Theorem}
\newtheorem{lem}{Lemma}[section]
\newtheorem{prop}{Proposition}[section]
\newtheorem{cor}{Corollary}[section]
\newtheorem{rem}{Remark}[section]

\begin{document}
\author{E.Feigin}
\title
[Bosonic formulas for affine branching functions]
{Bosonic formulas for affine branching functions}

\address{Evgeny Feigin:\newline
{\it Tamm Theory Division, Lebedev Physics Institute,
Russian Academy of Sciences,\newline 
Russia, 119991, Moscow, Leninski prospect, 53}\newline 
and \newline
{\it Independent University of Moscow,\newline 
Russia, Moscow, 119002, Bol'shoi Vlas'evskii, 11}}
\email{evgfeig@mccme.ru}

\begin{abstract}
In this paper we derive two bosonic (alternating sign) formulas for branching 
functions for general affine Kac-Moody Lie algebra $\g$. Both formulas are given 
in terms of Weyl group and string functions of $\g$.
\end{abstract}
\maketitle

\section*{Introduction}
Let $\g$ be an affine Kac-Moody Lie algebra, $\g_{fin}$ be the 
corresponding simple finite-dimensional algebra, $$\g=\g_{fin}\T\C[t,t^{-1}]
\oplus\C K\oplus\C d,$$ where $K$ is a central element and 
$[d,x\T t^i]=-ix\T t^i$ for $x\in\g.$ 
We fix the Cartan decomposition $\g=\n\oplus\h\oplus\n_-,$ with
$$\h=\h_{fin}\oplus\C K\oplus\C d,$$ where 
$\h_{fin}$ is the Cartan subalgebra of $\g_{fin}$.
We also denote $$\g'=[\g,\g]=\g\T\C[t,t^{-1}]\oplus\C K, \ \
\h'=\h_{fin}\oplus\C K.$$

Let $P^+_k\hk \h^*$ be the set of all dominant integrable level $k$ weights
of $\g$. We denote by $P^{'+}_k\hk {\h'}^*$ the image of $P_k^+$ with respect
to the restriction map $\h^*\to {\h'}^*$, $\la\mapsto\la'$. For $\la\in P_k^+$
we denote by $L_\la$ the corresponding irreducible highest weight $\g$ module 
and by $L_{\la'}$  a $\g'$ module which coincides with $L_\la$ as a vector
space and the action of $\g'$ is a restriction of the action of $\g$.

For $\la_1\in P_{k_1}^{'+}, \la_2\in P_{k_2}^{'+}$ consider a decomposition of 
the tensor product of $\g'$-modules
\begin{equation}
\label{LL}
L_{\la_1'}\T L_{\la_2'}=\bigoplus_{\mu'\in P_{k_1+k_2}^{'+}} 
L_{\mu'}\T C_{\la_1'\la_2'}^{\mu'}.
\end{equation}
Note that $C_{\la_1'\la_2'}^{\mu'}$ can be considered 
as a space of highest weight vectors of $\h'$-weight $\mu'$ in 
$L_{\la_1'}\T L_{\la_2'}$. 
To define a character of $C_{\la_1'\la_2'}^{\mu'}$ we assume 
$\la_1(d)=\la_2(d)=\mu(d)=0$ (note that for any $\la'\in P_k^{'+}$ there exists
$\bar\la\in P_k^+$ such that $\bar\la(d)=0$ and $\bar\la|_{\h'}=\la'$).
Then we obtain a grading by the operator $d$ on $L_{\la_1'}\T L_{\la_2'}$. Set
$$c_{\la_1'\la_2'}^{\mu'}(q)=\ch_q C_{\la_1'\la_2'}^{\mu'}=
\Tr\ q^d|_{C_{\la_1'\la_2'}^{\mu'}}.$$ 
These functions are called  $\g$ branching functions.
We note that in the conformal field theory branching functions appear as
characters of spaces of states of coset theories. These  characters differs
from $c_{\la_1'\la_2'}^{\mu'}(q)$ by an extra factor 
$q^{\triangle_{\la_1'}+\triangle_{\la_2'}-\triangle_{\mu'}}$, where
$\triangle_{\la'}$ is a conformal weight of $\la'$ (see \cite{DMS}).

There exist  different approaches to the study of $c_{\la_1'\la_2'}^{\mu'}(q)$
(see for example \cite{BNY, KMQ, R, DJKMO, S1, S2, SS, F, FOW}). 
These approaches give
different types formulas for some particular cases of branching functions.
In our paper we use homological technique to derive two bosonic formulas for 
$c_{\la_1'\la_2'}^{\mu'}(q)$ for general affine Kac-Moody algebras
(note that similar approach is utilized in \cite{FFJMT, FF, F}).
Let us briefly describe our results. Recall the Garland-Lepowsky theorem:
\begin{equation}
\label{G-L}
H_p(\n_-,L_{\mu})\simeq\bigoplus_{\substack{w\in W\\ l(w)=p}}\C_{w*\mu}, 
\end{equation}
where $W$ is a Weyl group of $\g$, $w*\mu$ is a shifted action of W and $l(w)$ 
is the length of $w$.
Note that $(\ref{G-L})$ is an isomorphism of $\h$-modules and $\C_{w*\mu}$ is
one-dimensional $\h$-module of the weight $w*\mu$.
From (\ref{G-L}) we obtain that for $\mu,\nu\in P_k^+$ homology
$H_p(\n_-,L_{\mu})^\nu$ (superscript denotes the corresponding $\h$-weight
subspace) vanishes if $p>0$ or $\mu\ne\nu$. In addition
$$H_0(\n_-,L_{\mu})^\mu\simeq\C_{\mu}.$$
Therefore from (\ref{LL}) we obtain
$$H_p(\n_-,L_{\la_1}\T L_{\la_2})^{\mu'}\simeq C_{\la_1'\la_2'}^{\mu'}
\delta_{0,p}$$
and  so
\begin{equation}
\label{mu'}
\sum_{p\ge 0}(-1)^p \ch_q H_p(\n_-,L_{\la_1}\T 
L_{\la_2})^{\mu'}= c_{\la_1'\la_2'}^{\mu'}(q).
\end{equation}
We now compute the same Euler characteristics using the BGG-resolution of
$L_{\la_1}$:
\begin{equation}
\ldots\to F_p\to\ldots\to F_0\to L_{\la_1}\to 0,
\end{equation}
where $F_p=\bigoplus_{l(w)=p} M_{w*\la_1}$ and 
$M_{w*\la_1}$ is the corresponding Verma module.
Tensoring the BGG-resolution by $L_{\la_2}$ we obtain the $\U(\n_-)$-free
resolution of $L_{\la_1}\T L_{\la_2}$:
\begin{equation}
\label{GG}
\ldots\to F_p\T L_{\la_2}\to\ldots\to F_0\T L_{\la_2}\to 
L_{\la_1}\T L_{\la_2}\to 0.
\end{equation}
Then the homology $H_p(\n_-,L_{\la_1}\T L_{\la_2})^{\mu'}$ can be counted as 
homology of a complex
\begin{equation}
\ldots\to \left[\C\T_{\U(\n_-)} (F_p\T L_{\la_2})\right]^{\mu'}\to\ldots\to 
\left[\C\T_{\U(\n_-)} (F_0\T L_{\la_2})\right]^{\mu'}\to 0.
\end{equation}
Therefore the Euler characteristics $(\ref{mu'})$ is given by the formula
\begin{equation}
\label{fhf}
\sum_{p\ge 0}(-1)^p \sum_{l(w)=p}q^{(w*\la_1)d}\ch_q 
(L_{\la_2})^{(\mu-w*\la_1)'}.
\end{equation}
We thus obtain our first bosonic formula for $c_{\la_1'\la_2'}^{\mu'}(q)$. 

To get the second formula we replace the "product" $L_{\la_1}\T
L_{\la_2}$ by the "fraction" $L_{\mu}\T L_{\la_1}^*$ and consider the
homology $H_p(\n_-,L_{\mu}\T L_{\la_1}^*)$, $\la_1\in P_{k_1}^+$,
$\mu\in P_{k_1+k_2}^+$.
We prove that
$$H_p(\n_-,L_{\mu}\T L_{\la_1}^*)^{\la_2'}=0 \text { for } p>0$$
and
$$H_0(\n_-,L_{\mu}\T L_{\la_1}^*)^{\la_2'}\simeq (C_{\la_1'\la_2'}^{\mu'})^*.$$
We thus obtain that
$$\sum_{p\ge 0} (-1)^p \ch_q H_p(\n_-,L_{\mu}\T L_{\la_1}^*)^{\la_2'}=
c_{\la_1'\la_2'}^{\mu'}(q^{-1}).$$
Using the BGG-resolution of $L_{\mu}$ we again rewrite this Euler 
characteristics in terms of the characters of spaces 
$(\C\T_{\U(\n_-)}(M_{w*\mu}\T L_{\la_1}^*))^{\la_2'}.$
This gives the following formula:
\begin{equation}
\label{shf}
c_{\la_1'\la_2'}^{\mu'}(q)=\sum_{p\ge 0}(-1)^p \sum_{l(w)=p}q^{-(w*\mu)(d)}
\ch_q (L_{\la_1})^{(w*\mu-\la_2)'}
\end{equation}
(recall that we assume $\la_1(d)=\la_2(d)=\mu(d)=0$). 
The specialization of this formula to the simplest case 
$\g=\slth$ gives the formula from \cite{BNY, KMQ, R} in the form of \cite{F}.
We note that $(\ref{shf})$ looks like  $(\ref{fhf})$, but the proof is 
much more complicated.

Our paper is organized as follows.

In Section $1$ we fix  affine Kac-Moody Lie algebras notations.

In Section $2$ we derive our first formula for branching functions 
$c_{\la_1'\la_2'}^{\mu'}(q)$ using the homology
$H_p(\n_-,L_{\la_1}\T L_{\la_2})^{\mu'}.$

In Section $3$ we derive our second formula for $c_{\la_1'\la_2'}^{\mu'}(q)$
using the homology of the "fraction" 
$H_p(\n_-,L_{\mu}\T L_{\la_1}^*)^{\la_2'}.$

In Section $4$ we specialize formulas from Sections $2$ and $3$ to the 
simplest case $\g=\slth$.

{\it Acknowledgements.}\quad This work was partially supported by
RFBR Grant 06-01-00037 and LSS 4401.2006.2.

\section{Affine Kac-Moody Lie algebras}
In this section we fix our notations on the affine Kac-Moody Lie algebras.
The main references are \cite{Kac, Kum}. 

Let $\g_{fin}$ be a simple finite-dimensional Lie algebra with the Cartan 
decomposition $\g_{fin}=\n_{fin}\oplus\h_{fin}\oplus (\n_-)_{fin}.$

Consider the corresponding affine algebra 
$$\g=\g_{fin}\T\C[t,t^{-1}]\oplus\C K\oplus\C d,$$
where $K$ is a central element and $[d,x\T t^i]=-ix\T t^i.$

We fix the Cartan decomposition $\g=\n\oplus\h\oplus\n_-,$
where 
\begin{gather*}
\n=\n_{fin}\T 1\oplus\g_{fin}\T t\C[t],\\
\h=\h_{fin}\oplus\C K\oplus\C d, \\
\n_-=(\n_-)_{fin}\T 1\oplus\g_{fin}\T t^{-1}\C[t^{-1}]
\end{gather*}
and denote $\g'=[\g,\g]=\g_{fin}\T\C[t,t^{-1}]\oplus\C K,$ 
$\h'=\h_{fin}\oplus\C K\hk\g'$.

Let $\al_i^{\vee}\in \h,\ \al_i\in\h^*,\ i=1,\ldots,n$,
be simple coroots and roots. Note that $\al_i^{\vee}$ form a basis of 
$\h_{fin}\oplus\C K.$
We denote by $\slt^{(i)}$ the $\slt$ Lie algebra spanned by $e_i,
\al_i^{\vee},f_i$, where $e_i,f_i,i=1,\ldots,n$ are the Chevalley
generators, $e_i\in\n, f_i\in\n_-.$
We note that 
$$\n=\bigoplus_{\al\in \triangle_+}\g_{\al},\ 
\n_-=\bigoplus_{\al\in \triangle_-}\g_{\al},$$
where 
$\triangle_+$ and $\triangle_-$ are the sets of positive and negative roots
and $\g_{\al}=\{x\in\g:[h,x]=\al(h)x\ \forall h\in\h\}.$
Spaces $\g_{\al_i}$ and $\g_{-\al_i}$ are spanned by $e_i$ and $f_i$.
Let 
\begin{equation}
\label{u-}
\un_-^{(i)}=\bigoplus_{\substack{\al\in \triangle_-\\ \al\ne -\al_i}}\g_{\al}.
\end{equation}
Note that 
$$\g_{-\al_i}\simeq\n_-/\un_-^{(i)}.$$

Let $P_k^+$ be the set of level $k$ integrable dominant $\g$-weights, i.e.
$$P_k^+=\{\la\in {\h}^*: \la(\al_i^{\vee})\in\Z_{\ge 0},\ \la(K)=k\}.$$
We also denote by $P_k^{'+}\hk {\h'}^*$ the image of $P_k^+$ with respect to
the projection $\h^*\to {\h'}^*$, $\la\mapsto\la'$. 
For $\la\in P_k^+$ let $L_{\la}$ be an integrable
highest weight $\g$-module with highest weight vector $v_{\la}\in L_{\la}$
such that 
$$\n v_{\la}=0,\ \U(\n_-)v_{\la}=L_{\la},\ h(v_{\la})=\la(h)v_{\la}, h\in\h.$$
Let $L_{\la'}$ be $\g'$ module which coincides with $L_\la$ as a vector space 
and the action of $\g'$ is a restriction of the action of $\g$.
For any $\al\in {\h'}^*$ set
$$(L_{\la})^{\al}=\{v\in L_\la:\ hv=\al(h)v\ \forall h\in\h'\}.$$ 
Note that $L_\la$ is graded by an operator $d$. We set 
$$\ch_q (L_{\la})^{\al}=\Tr q^d|_{(L_{\la})^{\al}}.$$

Fix  $\la_1'\in P_{k_1}^{'+},\ \la_2'\in P_{k_2}^{'+},$
and consider the decomposition of the tensor product of $\g'$-modules:
\begin{equation}
\label{tp}
L_{\la_1'}\T L_{\la_2'}=
\bigoplus_{\mu'\in P_{k_1+k_2}^{'+}} 
C_{\la_1'\la_2'}^{\mu'} \T L_{\mu'}.
\end{equation}
The space  $C_{\la'_1\la'_2}^{\mu'}$ 
can be identified with a subspace of highest weight vectors of 
$\h'$-weight $\mu'$ in $L_{\la_1'}\T L_{\la_2'}.$
To define a character of $C_{\la'_1\la'_2}^{\mu'}$  one needs to
fix an action of the operator $d$ on each $\g'$ module $L_{\la'}$.
Note that if $\la|_{\h'}=\bar\la|_{\h'}$ then $L_{\la'}\simeq L_{\bar\la'}$.
Therefore an action of $d$ on $L_{\la'}$ depends on the choice of $\la(d)$. 
It is convenient for us to choose a normalization 
$\la_1(d)=\la_2(d)=\mu(d)=0$.
This defines the characters of $L_{\la_1'}\T L_{\la_2'}$ and of 
$C_{\la'_1\la'_2}^{\mu'}$. The character 
$\ch_q C_{\la'_1\la'_2}^{\mu'}$ is called $\g$ branching function and 
is denoted by $c_{\la_1'\la_2'}^{\mu'}(q)$:
$$c_{\la_1'\la_2'}^{\mu'}(q)=\Tr q^d|_{C_{\la_1'\la_2'}^{\mu'}}.$$

Recall that spaces $C_{\la_1' \la_2'}^{\mu'}$ appear in the conformal field theory
as spaces of states of coset theories (see \cite{DMS}). 
Namely the Sugawara construction defines an action of the Virasoro algebra
with generators $L_n$ on each $L_{\la'}$.
In particular for the operator $L_0$ one has
$$L_0 v_{\la}=\triangle_{\la'} v_{\la},\ [L_0,x\T t^i]=-ix\T t^i$$ 
($\triangle_{\la'}$ is a conformal weight).
Now the GKO construction (see \cite{GKO}) defines an action of Vir on the 
tensor product $L_{\la'_1}\T L_{\la'_2}$ which commutes with the diagonal
action of $\g'$.
Namely one puts
$$L_n^{GKO}=L_n^{(1)}\T \Id +\Id\T L_n^{(2)}-L_n^{diag},$$
where $L_n^{(1)},\ L_n^{(2)}$ and $L_n^{diag}$ are Sugawara operators acting
on $L_{\la_1'},\ L_{\la_2'}$ and $L_{\la_1'}\T L_{\la_2'}$ respectively.
Therefore, we obtain a structure of Vir-module on
$C_{\la_1'\la_2'}^{\mu'}$ and an equality
$$\Tr q^{L_0}|_{C_{\la_1'\la_2'}^{\mu'}}=c_{\la_1'\la_2'}^{\mu'}(q)
q^{\triangle_{\la_1'}+\triangle_{\la_2'}-\triangle_{\mu'}},$$
where the left hand side is a character of the space of states of the 
corresponding coset model.

In the end of this section we recall the Weyl group notations, 
the Garland-Lepowsky $\n_-$-homology
theorem and the BGG resolution.
Let $W$ be the Weyl group of $\g$, generated by simple reflections $s_i$.
We denote by $l(w)$ the length of an element $w\in W$.
Recall that the shifted action of $W$ on $\h^*$ is given by
$w*\la =w(\la+\rho)-\rho$, where $\rho(\al_i^{\vee})=1.$ 
We will need the following lemma: 
\begin{lem}
\label{Wlem}
$a)$\ If $\la\in P^+_k$ and $w*\la\in P^+_k$ then $w=e$. \\
$b)$\ If $(w*\la)\al_i^\vee\le -1$ for some $\la\in P^+_k$ then  
$l(s_iw)<l(w)$.
\end{lem}
\begin{proof}
To prove $a)$ we rewrite an equality $w*\la=\la_1$ as 
$w(\la+\rho)=\la_1+\rho$. But
if $\la,\la_1\in P^+_k$ then 
$$(\la+\rho)\al_i^\vee>0, \quad (\la_1+\rho)\al_i^\vee>0\quad \forall i.$$
Therefore both  $\la+\rho$ and $\la_1+\rho$ are the elements of the Weyl dominant 
chamber and so $\la=\la_1$. But from $w*\la=\la$ one gets $w=e$ (see 
Lemma $3.2.5$ from \cite{Kum}).

We now prove $b)$. Note that  $(w*\la)\al_i^\vee\le -1$ is equivalent to
$$(\la+\rho)(w^{-1} \al_i^\vee)\le 0.$$ This gives 
$w^{-1} \al_i^\vee=\sum_{i=1}^n c_j  \al_j^\vee$ with $c_j\le 0$. Therefore
$l(w^{-1}s_i)<l(w)$. Lemma is proved. 
\end{proof}

The following theorem is proved in \cite{GL}: 
\begin{theo*}
For any $\la\in P_k^+$ we have an isomorphism of $\h$-modules:
$$H_p(\n_-,L_\la)\simeq\bigoplus_{\substack{w\in W\\ l(w)=p}}\C_{w*\la},$$
where $\C_{w*\la}$ is one-dimensional $\h$-module of the weight $w*\la$.
\end{theo*}
We will also need the BGG resolution of integrable irreducible representations
$L_{\la}$ (see \cite{BGG, Kum}). 
Namely there exists an exact sequence of $\g$-modules and $\g$-homomorphism
\begin{equation}
\label{BGGres}
\ldots\to\bigoplus_{\substack{w\in W\\l(w)=p}}M_{w*\la}
\to\ldots\to M_{\la}\to L_{\la}\to 0,
\end{equation}
where $M_{\mu}$ is the weight $\mu$ Verma module.

\section{First homological bosonic formula}

\begin{lem}
For any $\la_1\in P^+_{k_1}$, $\la_2\in P^+_{k_2}$ and 
$\mu\in P^+_{k_1+k_2}$ 
we have:
\begin{gather*} 
C_{\la_1'\la_2'}^{\mu'}\simeq H_0(\n_-, L_{\la_1}\T L_{\la_2})^{\mu'},\\
H_p(\n_-, L_{\la_1}\T L_{\la_2})^{\mu'}=0 \text{ for all } p>0.
\end{gather*}
\end{lem}
\begin{proof}
Recall that for any $\la\in P^+_k$ and $e\ne w\in W$  one has
$w*\la\notin P^+_k$ and so $(w*\la)'\notin P^{'+}_k$. Therefore, from 
Garland-Lepowsky theorem we obtain  that $H_p(\n_-, L_\la)^\mu=0$
unless $p=0$ and $\la=\mu$. Now our lemma follows from the decomposition 
$(\ref{tp})$ and an equality 
$$H_0(\n_-, L_{\la_1}\T L_{\la_2})^{\mu'}\simeq
\bigoplus_{\substack{\bar\mu\in P^+_k\\ \bar\mu|_{\h'}=\mu'}} 
H_0(\n_-, L_{\la_1}\T L_{\la_2})^{\bar\mu}.$$
\end{proof}

\begin{cor}
\label{br}
For any $\la_1\in P^+_{k_1}$, $\la_2\in P^+_{k_2}$ and $\mu\in P^+_{k_1+k_2}$
one has
\begin{equation}
\label{ec}
c_{\la_1'\la_2'}^{\mu'}(q)= 
\sum_{p\ge 0} (-1)^p \ch_q H_p(\n_-, L_{\la_1}\T L_{\la_2})^{\mu'}.
\end{equation}
\end{cor}

We now compute  the Euler characteristics $(\ref{ec})$ using the 
BGG-resolution of $L_{\la_1}$. 
Tensoring $(\ref{BGGres})$ by $L_{\la_2}$ we obtain the $\U(\n_-)$-free 
resolution of 
$L_{\la_1}\T L_{\la_2}$. Therefore the following complex counts 
$H_p(\n_-, L_{\la_1}\T L_{\la_2})$:
\begin{equation}
\label{coinv}
\ldots \to \C\T_{\U(\n_-)} (L_{\la_2}\T F_1) \to 
\C\T_{\U(\n_-)} (L_{\la_2}\T F_0)\to 0,
\end{equation}
where $F_p=\bigoplus_{l(w)=p} M_{w*\la_1}$. 
We can rewrite 
($\ref{coinv}$) as
\begin{equation}
\label{cu}
\ldots \to (\C\T_{\U(\n_-)} F_1)\T L_{\la_2} \to 
(\C\T_{\U(\n_-)} F_0)\T L_{\la_2}\to 0.
\end{equation}

\begin{lem}
\label{form}
\begin{multline*}
\sum_{p\ge 0} (-1)^p \ch_q H_p(\n_-, L_{\la_1}\T L_{\la_2})^{\mu'}=\\
\sum_{p\ge 0} (-1)^p \sum_{l(w)=p} q^{(w*\la_1)d} 
\ch_q (L_{\la_2})^{(\mu-w*\la_1)'}.
\end{multline*}
\end{lem}
\begin{proof}
Recall that $F_p=\bigoplus_{l(w)=p} M_{w*\la_1}$. Therefore,
$$\C\T_{\U(\n_-)} F_p= \bigoplus_{l(w)=p} \C_{w*\la_1}.$$
Now our lemma follows from the equality of Euler characteristics of the complex
$(\ref{cu})$ and the right hand side of $(\ref{ec})$. 
\end{proof}

\begin{prop}
\label{mainhom}
We have a bosonic formula for the branching functions:
\begin{equation}
\label{genhom}
c_{\la_1'\la_2'}^{\mu'}(q)=
   \sum_{p\ge 0} (-1)^p \sum_{\substack{w\in W\\ l(w)=p}} q^{(w*\la_1)d} 
\ch_q (L_{\la_2})^{(\mu-w*\la_1)'}
\end{equation}
\end{prop}
\begin{proof}
Follows from Corollary $\ref{br}$ and Lemma $\ref{form}$. 
\end{proof}

\begin{rem}
\label{rem}
We can use the BGG resolution of $L_{\la_2}$ instead of 
$L_{\la_1}$. This interchanges $\la_1$ and $\la_2$ in the right hand 
side of $(\ref{genhom})$ and leads to another formula for branching 
functions 
$c_{\la_1'\la_2'}^{\mu'}(q).$ 
\end{rem}

\section{Second homological bosonic formula}
In this section we study homology $H_p(\n_-, L_\mu\T L_{\la_1}^*)$ replacing 
the tensor product  $L_{\la_1}\T L_{\la_2}$ from the previous section by the 
"fraction" $L_\mu\T L_{\la_1}^*$. We note that though  $L_\mu\T L_{\la_1}^*$
does not belong to the category $\EuScript{O}$ (the eigenvalues of the 
operator $d$ are not bounded from below) it is still integrable. So we first
prove some statements about integrable representations.

Recall the definition $(\ref{u-})$ of the subalgebra $\un^{(i)}_-$.
\begin{lem}
\label{nu}
Let $M$ be an integrable $\g$ module. Then
$$H_n(\n_-, M)\simeq H_0(\g_{-\al_i}, H_n(\un^{(i)}_-, M))\oplus 
H_1(\g_{-\al_i}, H_{n-1}(\un^{(i)}_-,M)).$$
\end{lem}
\begin{proof}
We consider the Hochschild-Serre spectral sequence associated with a pair
$\un^{(i)}_-\hk\n_-$. 
Note that $\un^{(i)}_-$ is an ideal and 
$\n_-/\un^{(i)}_-\simeq \g_{-\al_i}$. The second
term of this spectral sequence is given by 
$$E^2_{p,q}=H_p(\g_{-\al_i}, H_q(\un^{(i)}_-,M)).$$
We prove our lemma by showing that $E^2_{p,q}=E^{\infty}_{p,q}$. 

Because of the integrability condition $M$ is a direct sum of irreducible 
finite-dimensional $\slt^{(i)}$ modules. Therefore the same is true for
$\Lambda^{q}(\un^{(i)}_-)\T M$ and also for $H_q(\un^{(i)}_-,M)$. 
For any nonnegative integer $s$ we denote by $\pi_s$ an irreducible
$\slt$ module with highest weight $s$ ($\dim\pi_s=s+1$) and fix 
highest and lowest weight vectors $v_s$ and $u_s$.
Let 
$\pi_s\hk H_q(\un^{(i)}_-,M)$ be a direct summand and
$$\alpha_p\in \Lambda^p(\g_{-\al_i})\T \pi_s$$ be a chain 
representing  some class in
$H_p(\g_{-\al_i},\pi_s)\hk H_p(\g_{-\al_i}, H_q(\un^{(i)}_-,M)).$
We set
$$\al_0=v_s, \ \ \ \al_1=f_i\T u_s.$$
Let 
$\be_p\in \Lambda^p(\g_{-\al_i})\T \Lambda^{q}(\un^{(i)}_-)\T M$  
be the chains of the form
$$\be_0=x_0,\ \ \ \be_1=f_i\T x_1$$
which represent $\al_p$ (i.e. $x_0$ represents $v_s$ and $x_1$ represents $u_s$).
Now let $d_{\n_-}$ ($d_{\un^{(i)}_-}$) be the differential in the standard 
complex for $H_n(\n_-,M)$ ($H_n(\un^{(i)}_-,M)$). 
We state that $d_{\n_-} \be_p =0$. In fact, for 
$p=0$ this just follows from $v_s\in H_q(\un^{(i)}_-, M)$. Now let
$p=1$. Then 
$$d_{\n_-} \be_1= d_{\n_-} (f_i\T x_1).$$
We know that $d_{\un^{(i)}_-} x_1=0$ and $f_ix_1=0$ (because $x_1$ represents
the lowest weight vector). This gives $d_{\n_-} \be_1=0$.
But because of
$d_{\n_-}\be=0$ we obtain that differentials $d_2,d_3,\ldots$ in the 
Hochschild-Serre spectral sequence
are trivial and $E^2_{p,q}=E^\infty_{p,q}$. Lemma is proved.  
\end{proof}

\begin{cor}
\label{H0w}
Let $M$ be an integrable level $k$ $\g$ module. Then 
$$H_0(\n_-,M)^\la=0 \text{ unless } \la\in P^+_k.$$
\end{cor}
\begin{proof}
Because of Lemma \ref{nu} we obtain
$$H_0(\n_-,M)\simeq H_0(\g_{-\al_i}, H_0(\un^{(i)}_-,M))$$
for all $i=1,\ldots,n$. We recall that $H_0(\g_{-\al_i},\pi_s)$ is
one-dimensional space of the $\al_i^\vee$-weight $s$. Therefore, because
$H_0(\un^{(i)}_-,M)$ is a direct sum of finite-dimensional $\slt^{(i)}$ modules, we 
obtain $\lambda(\al_i^\vee)\in\Z_{\ge 0}$ for any weight $\lambda$ of 
$H_0(\n,M)$. This gives $\la\in P_k^+$.
\end{proof}

\begin{prop}
\label{wla}
Let $M$ be an integrable level $k$ $\g$ module. Then\\
$a)$ $H_n(\n_-,M)^{w*\la}=0$ for $w\in W$, $\la\in P^+_k$ if $l(w)>n$.\\
$b)$ $H_n(\n_-,M)^{w*\la}\simeq H_{n-l(w)}(\n_-,M)^\la$ for any 
$w\in W$, $\al\in P^+_k$ such that $l(w)\le n$.
\end{prop}
\begin{proof}
We prove $a)$ and $b)$  simultaneously using  
$$H_n(\n_-,M)^\mu\simeq H_0(\g_{-\al_i}, H_n(\un^{(i)}_-,M))^\mu\oplus 
H_1(\g_{-\al_i}, H_{n-1}(\un^{(i)}_-,M))^\mu$$ 
and the induction on $n$. The case $n=0$
follows from Corollary \ref{H0w}. Suppose our lemma is proved for $m<n$. We
assume $l(w)>0$ (otherwise $a)$ and $b)$ are trivial). Then there exists
$i$ such that 
\begin{equation}
\label{notin}
(w*\la)\al_i^\vee \le -1.
\end{equation}
We have
$$H_n(\n_-,M)^{w*\la}\simeq H_0(\g_{-\al_i}, H_n(\un^{(i)}_-,M))^{w*\la}\oplus   
H_1(\g_{-\al_i}, H_{n-1}(\un^{(i)}_-,M))^{w*\la}.$$
Because of the condition $(\ref{notin})$ we have 
$$H_0(\g_{-\al_i}, H_n(\un^{(i)}_-,M))^{w*\la}=0$$
and therefore
\begin{multline}
\label{n}
H_n(\n,M)^{w*\la}\simeq H_1(\g_{-\al_i}, H_{n-1}(\un^{(i)}_-,M))^{w*\la}\simeq\\
H_0(\g_{-\al_i}, H_{n-1}(\un^{(i)}_-,M))^{(s_iw)*\la},
\end{multline}
because 
$$H_1(\g_{-\al_i},\pi_s)^\al\simeq H_0(\g_{-\al_i},\pi_s)^{s_i*\al}$$
for any $\al$.
We also know that
\begin{multline*}
H_{n-1}(\n_-,M)^{(s_iw)*\la}\simeq 
H_0(\g_{-\al_i}, H_{n-1}(\un^{(i)}_-,M))^{(s_iw)*\la}\oplus\\   
H_1(\g_{-\al_i}, H_{n-2}(\un^{(i)}_-,M))^{(s_iw)*\la}
\end{multline*}
and because of $((s_iw)*\la)\al_i^\vee=-(w*\la)\al_i^\vee-2\ge -1$ 
(see ($\ref{notin}$))
we obtain
$$H_1(\g_{-\al_i}, H_{n-2}(\un,M))^{(s_iw)*\la}=0$$
(because $H_1(\g_{-\al_i}, \pi_s)^t=0$ for any $t\ge -1$).
Therefore 
\begin{equation}
\label{n-1}
H_{n-1}(\n_-,M)^{(s_iw)*\la}\simeq 
H_0(\g_{-\al_i}, H_{n-1}(\un^{(i)}_-,M))^{(s_iw)*\la}.
\end{equation}
From  $(\ref{n})$ and $(\ref{n-1})$ we obtain
\begin{equation}
\label{reduce}
H_n(\n_-,M)^{w*\la}\simeq H_{n-1}(\n_-,M)^{(s_iw)*\la}.
\end{equation}
Note that because of $(\ref{notin})$ and Lemma $\ref{Wlem}$  
$l(s_iw)=l(w)-1$. 

Now suppose that $n<l(w)$. Then iterating 
$(\ref{reduce})$ we obtain
$$H_n(\n_-,M)^{w*\la}\simeq H_0(\n_-,M)^{w'*\la}$$
for some $w'$ with $l(w')>0$. But this homology vanishes because of 
Corollary $\ref{H0w}$. This gives $a)$. To obtain $b)$ one needs to iterate 
$(\ref{reduce})$. Proposition is proved.
\end{proof}

Let $\omega:\g\to\g$ be the Chevalley involution defined by $e_i\to -f_i$, 
$f_i\to -e_i$, $h\to -h$ ($h\in\h$). 
For $\g$ module $V$ we denote by $V^\omega$ a $\g$ module which coincides
with $V$ as a vector space and the action of $\g$ is twisted by $\omega$.

\begin{lem}
\label{ng}
Let $M$ be some $\g$ module, $\la\in\h^*$. Then
$$H_n(\n_-, M)^\la\simeq  H_n(\g,\h, M\T M_\la^\omega),$$
where $M_\la$ is the Verma module.
\end{lem}
\begin{proof}
We first rewrite
$$H_n(\n_-, M)^\la\simeq H_n(\bo_-,\h, M\T \C_{-\la}),$$
where $\C_{-\la}$ is one-dimensional $\bo_-=\n_-\oplus\h$-module with trivial 
action of $\n_-$. Now our lemma follows from 
$$\mathrm{Ind}_{\bo_-}^{\g} (M\T \C_{-\la})\simeq M\T M_\la^\omega.$$ 
\end{proof}

In what follows we study homology  
$H_n(\n_-, L_\mu\T L_{\la_1}^*)^{\la_2}$ for the triple of weights
$\la_1\in P^+_{k_1}$, $\la_2\in P^+_{k_2}$, $\mu\in P^+_{k_1+k_2}$. 
Note that $L_\mu\T L_{\la_1}^*$ is integrable $\g$ module.
Because of Lemma $\ref{ng}$ we have an isomorphism
$$H_n(\n_-, L_\mu\T L_{\la_1}^*)^{\la_2}\simeq 
H_n(\g,\h, L_\mu\T L_{\la_1}^* \T M_{\la_2}^\omega).$$

\begin{prop}
\label{barE}
There exists a spectral sequence $\bar E^{r}_{p,q}$ with 
$$\bar E^1_{p,q}=
\bigoplus_{w:\ l(w)=p} H_q(\n_-, L_\mu\T L_{\la_1}^*)^{w*\la_2},$$
such that 
$\bar E^r_{p,q}$ converges to 
$H_\bullet(\g,\h, L_\mu\T L_{\la_1}^* \T L_{\la_2}^*).$
In addition $\bar E^1_{p,q}=0$ for $p>q$.
\end{prop}
\begin{proof}
We first note that $\bar E^1_{p,q}=0$ for $p>q$ because of part $a)$ of 
Proposition $\ref{wla}$. 

Now
consider the BGG resolution
$$\ldots \to F_{\la}(2)\to F_{\la}(1)\to F_{\la}(0)\to L_{\la}\to 0,\
F_{\la}(p)=\bigoplus_{w:\ l(w)=p} M_{w*\la}.$$
Recall that for any $\la\in P^+_k$ $(L^*_{\la})^\omega\simeq L_\la$.
We thus obtain the dual BGG-resolution
\begin{equation*}
\ldots \to F^\omega_{\la}(2)\to F^\omega_{\la}(1)\to 
F^\omega_{\la}(0)\to L^*_{\la}\to 0,\
F^\omega_{\la}(p)=\bigoplus_{w:\ l(w)=p} M^\omega_{w*\la}.
\end{equation*}
This gives the following resolution 
\begin{equation}
\label{dBGG}
\ldots \to L_\mu\T L_{\la_1}^*\T F^\omega_{\la_2}(1)\to 
L_\mu\T L_{\la_1}^*\T F^\omega_{\la_2}(0)\to 
L_\mu\T L_{\la_1}^*\T L^*_{\la_2}\to 0.
\end{equation}

In order to establish a connection between 
$$H_q(\g,\h, L_\mu\T L_{\la_1}^* \T L_{\la_2}^*)\text{ and }
H_q(\n_-, L_\mu\T L_{\la_1}^*)^{w*\la_2}$$
we use a bi-complex $K_{p,q}$ associated with the resolution $(\ref{dBGG})$:
$$K_{p,q}=
\left[ \Lambda^q(\g/\h)\T L_\mu\T L_{\la_1}^* \T F_{\la_2}^\omega(p)\right]
^0
$$
(the space of $\h$-invariants).
The first term of the corresponding spectral sequence is given by
$${\bar E}^1_{p,q}=
H_q(\g, \h, L_\mu\T L_{\la_1}^* \T F_{\la_2}^\omega(p))$$
and ${\bar E}^r_{p,q}$ converges to 
$H_n(\g,\h, L_\mu\T L_{\la_1}^* \T L_{\la_2}^*).$
From the definition of $F_{\la_2}^\omega(p)$ and Lemma $\ref{ng}$ we obtain
$$\bar E^1_{p,q}=
\bigoplus_{w:\ l(w)=p} H_q(\n_-, L_\mu\T L_{\la_1}^*)^{w*\la_2}.$$
Proposition is proved.
\end{proof}

\begin{cor}
\label{H0}
$H_0(\n_-, L_\mu\T L_{\la_1}^*)^{\la_2'}\simeq (C_{\la_1'\la_2'}^{\mu'})^*$.
\end{cor}
\begin{proof}
Note that $H_0(\n_-, L_\mu\T L_{\la_1}^*)^{\la_2'}$ is isomorphic to
$$\bigoplus_{\substack{\bar \la_2\in P_k^+\\ \bar \la_2|_{\h'}=\la_2'}}
H_0(\n_-, L_\mu\T L_{\la_1}^*)^{\bar \la_2}.$$
Fix some $\bar\la_2$ with $\bar \la_2|_{\h'}=\la_2'$.
Then using the spectral sequence from the Proposition $\ref{barE}$ we obtain 
$$H_0(\n_-, L_\mu\T L_{\la_1}^*)^{\bar\la_2}=\bar E^1_{0,0}$$
and $\bar E^1_{n,0}=0$ for all $n>0$. Therefore 
$\lim_{r\to\infty} \bar E^r_{0,0}=\bar E^1_{0,0}$.
In addition 
$$\bar E^\infty_{0,0}\simeq H_0(\g,\h, L_\mu\T L_{\la_1}^* \T L_{\bar\la_2}^*)
\simeq H^0(\g,\h, L_\mu^*\T L_{\la_1} \T L_{\bar \la_2})^*.$$
Now our Corollary follows from
$$H^0(\g,\h, L_\mu^*\T L_{\la_1} \T L_{\bar\la_2})^*\simeq 
\mathrm{Hom}_{\g} (L_\mu, L_{\la_1} \T L_{\bar\la_2})^*$$ and
$$\bigoplus_{\substack{\bar \la_2\in P_k^+\\ \bar \la_2|_{\h'}=\la_2'}}
H^0(\g,\h, L_\mu^*\T L_{\la_1} \T L_{\bar\la_2})^*\simeq 
\mathrm{Hom}_{\g'} (L_{\mu'}, L_{\la_1'} \T L_{{\bar\la_2}'})^*.$$
\end{proof}

We now study the special case of $H_n(\n_-, L_\mu\T L_{\la_1}^*)^{\la_2}$
with $\la_2=0$ and $\mu=\la_1$.
\begin{lem}
\label{mumu*}
$H_n(\n_-,L_\mu\T L_\mu^*)^0=0$ for $n>0$.
\end{lem}
\begin{proof}
We consider a filtration $(L_\mu^*)_m$ on $L_\mu^*$:
$$(L_\mu^*)_m=\mathrm{span}\bra e_{i_1}\ldots e_{i_s} v_\mu^*,\ s\le m, 
1\le i_l\le n \ket,$$
where $v_\mu^*$ is a lowest weight vector of $L_\mu^*$.
This induces a filtration 
$$(\left[\Lambda^n(\n_-)\T L_\mu\T L_\mu^*\right]^0)_m =
\left[\Lambda^n(\n_-)\T L_\mu\T (L_\mu^*)_m\right]^0.$$
For the associated spectral sequence one has 
$$E^1_{n,m}=\left[ H_{n+m}(\n_-,L_\mu)\T (L_\mu^*)_m/(L_\mu^*)_{m-1}\right]^0.$$
Because of 
$H_{n+m}(\n_-,L_\mu)\simeq\bigoplus_{w:\ l(w)=n+m} \C_{w*\mu}$ we obtain
$$E^1_{n,m}=\bigoplus_{w:\ l(w)=n+m}
\left[ (L_\mu^*)_m/(L_\mu^*)_{m-1}\right]^{-w*\mu}.$$
But 
$$(L_\mu^*)^{-w*\mu}\simeq (L_\mu^*)^{-w\mu+\rho-w\rho}\simeq
(L_\mu^*)^{-\mu+w^{-1}\rho-\rho}=0$$
because $w^{-1}\rho-\rho < 0$ for $l(w)>0$. This gives 
$E^1_{n,m}=0$ for $n+m\ne 0$. Lemma is proved.
\end{proof}

In the following Lemma we calculate 
homology $H_n(\g,\h, L_\mu\T L_\nu^*)$ for two weights
$\mu,\nu\in P^+_{k_1+k_2}$.

\begin{lem}
\label{munu}
Let $\mu,\nu\in P^+_{k_1+k_2}$. Then
$$\dim H_{2n}(\g,\h, L_\mu\T L_\nu^*)=\delta_{\mu,\nu} \#\{w\in W:\ l(w)=n\}.$$
In addition $H_{2n-1}(\g,\h, L_\mu\T L_\nu^*)=0$.
\end{lem}
\begin{proof}
Because of the isomorphism $(L_\nu^*)^\omega\simeq L_\nu$ the BGG resolution
gives the following resolution:
$$\ldots\to \bigoplus_{\substack{w_1,w_2\in W \\ l(w_1)+l(w_2)=p}} 
M_{w_1*\mu}\T M_{w_2*\nu}^\omega\to\ldots \to  
M_\mu\T M_\nu^\omega\to L_\mu\T L_\nu^*\to 0.
$$
This resolution is $(\g,\h)$-free. Therefore a complex 
$G_\bullet$ with
$$G_p= \bigoplus_{\substack{w_1,w_2\in W\\ l(w_1)+l(w_2)=p}} 
\left[\C\T_{U(\g/\h)} (M_{w_1*\mu}\T M_{w_2*\nu}^\omega) \right]^0$$
counts $H_p(\g,\h, L_\mu\T L_\nu^*)$.
Note that $G_p=\bigoplus_{\substack{l(w_1)+l(w_2)=p\\ w_1*\mu= w_2*\nu}} \C_0$.
In view of $\mu,\nu\in P^+_{k_1+k_2}$ the condition
$w_1*\mu=w_2* \nu$ is equivalent to $w_1=w_2$, $\nu=\mu$. Therefore
$G_{2n-1}=0$, $n\ge 1$ and 
$\dim G_{2n}= \delta_{\mu\nu} \#\{w\in W:\ l(w)=n\}.$
Lemma is proved.
\end{proof}

\begin{prop}
\label{triv}
The natural map
$$\tau: H_n(\g,\h, L_\mu\T L_{\la_1}^*\T M_{\la_2}^\omega) \to
        H_n(\g,\h, L_\mu\T L_{\la_1}^*\T L_{\la_2}^*)$$
is trivial for $n>0$. 
\end{prop}
\begin{proof}
Because of an isomorphism 
$$H_n(\g,\h, L_\mu\T L_{\la_1}^*\T M_{\la_2}^\omega)\simeq 
H_n(\bo_-,\h, L_\mu\T L_{\la_1}^*\T \C_{\la_2}^*)$$
it is enough to prove that the map 
$$\tau': H_n(\bo_-,\h, L_\mu\T L_{\la_1}^*\T \C_{\la_2}^*)\to
H_n(\g,\h, L_\mu\T L_{\la_1}^*\T L_{\la_2}^*)$$
is trivial.
Note that $\tau'$ comes from the natural embedding
\begin{equation}
\label{imath}
\imath: 
\left[ \Lambda^n(\bo_-/\h)\T L_\mu\T L_{\la_1}^*\T \C_{\la_2}^*\right]^0\hk
\left[ \Lambda^n(\g/\h)\T L_\mu\T L_{\la_1}^*\T L_{\la_2}^*\right]^0.
\end{equation}
Therefore it suffices to show that any chain 
$$c\in \left[ \Lambda^n(\bo_-/\h)\T L_\mu\T L_{\la_1}^*\T L_{\la_2}^*
\right]^0\hk 
\left[ \Lambda^n(\g/\h)\T L_\mu\T L_{\la_1}^*\T L_{\la_2}^*\right]^0$$
defines a trivial class in $H_n(\g,\h, L_\mu\T L_{\la_1}^*\T L_{\la_2}^*)$.

Because of Lemma $\ref{munu}$ it is enough to show that any chain from
$$\left[ \Lambda^n(\n_-)\T L_\mu\T L_\mu^*\right]^0\hk 
\left[ \Lambda^n(\g/\h)\T L_\mu\T L_{\la_1}^*\T L_{\la_2}^*
\right]^0$$
defines a trivial element in $H_n(\g,\h, L_\mu\T L_{\la_1}^*\T L_{\la_2}^*)$
for $n>0$. But from Lemma $\ref{mumu*}$ we know that
$$H_n(\n_-, L_\mu\T L_\mu^*)^0=0 \text{ for } n>0.$$
Proposition is proved.
\end{proof}

\begin{theo}
\label{main}
For any $\la_1\in P^+_{k_1}$, $\la_2\in P^+_{k_2}$, $\mu\in P^+_{k_1+k_2}$
and $n>0$
$$H_n(\n_-, L_\mu\T L_{\la_1}^*)^{\la_2}=0.$$ 
\end{theo}
\begin{proof}
We use the spectral sequence from Proposition $\ref{barE}$:
$$\bar E^1_{p,q}=\bigoplus_{w:\ l(w)=p} 
H_q(\n_-, L_\mu\T L_{\la_1}^*)^{w*\la_2}\simeq
\bigoplus_{l(w)=p} H_q(\g,\h, L_\mu\T L_{\la_1}^*\T M^\omega_{w*\la_2}).$$
From Corollary $\ref{H0}$ we obtain $\bar E^1_{0,0}\simeq 
(C_{\la'_1\la'_2}^{\mu'})^*$ and therefore 
\begin{equation}
\label{E1kk}
\bar E^1_{k,k}\simeq\bigoplus_{w:\ l(w)=k} (C_{\la'_1\la'_2}^{\mu'})^*.
\end{equation}
In addition Propositions $\ref{wla}$ and $\ref{barE}$ gives
\begin{equation}
\label{know}
\bar E^1_{p,n+p}\simeq \bar E^1_{0,n},\
\bar E^1_{p,q}=0 \text{ for } p>q.
\end{equation}
Our goal is to show that $\bar E^1_{0,n}=0$, $n>0$. 
(Because of $(\ref{know})$ this is
equivalent to the proof of $\bar E^1_{p,q}=0$ for $p\ne q$). 
Note that this agrees with a fact that $\bar E^1_{k,k}$ is isomorphic to
$H_{2k}(\g,\h, L_\mu\T L_{\la_1}^*\T L_{\la_1}^*)$ (see Lemma $\ref{munu}$ and 
$(\ref{E1kk})$).
  
We prove  the statement $\bar E^1_{0,2n-1}=0$, $\bar E^1_{0,2n}=0$ by
induction on $n\ge 1$. First let $n=1$. Because of Proposition $\ref{triv}$ 
the map $$\bar E^1_{0,2}\to H_2(\g,\h, L_\mu\T L_{\la_1}^*\T L_{\la_1}^*)$$
is trivial and therefore $\bar E^\infty_{0,2}=0$. This gives 
$\bar E^\infty_{1,1}=\bar E^1_{1,1}$, because 
$$\bar E^1_{1,1}\simeq H_2(\g,\h, L_\mu\T L_{\la_1}^*\T L_{\la_1}^*)\simeq
\bar E^\infty_{1,1}\oplus \bar E^\infty_{2,0}.$$
So a differential 
$d_1: \bar E^1_{1,1}\to \bar E^1_{0,1}$ is trivial and  
$$\bar E^1_{0,1}=\bar E^\infty_{0,1}= 
H_1(\g,\h, L_\mu\T L_{\la_1}^*\T L_{\la_1}^*)=0.$$
According to Proposition $\ref{wla}$ we obtain $\bar E^1_{1,2}=0$, which
gives 
$$\bar E^1_{0,2}= \bar E^\infty_{0,2}= 0.$$ 

Now suppose $\bar E^1_{0,s}=0$ for $s\le 2(n-1)$. This gives 
$\bar E^1_{p,s+p}=0$ for $s\le 2(n-1)$, $p\ge 0$. Note that the map
$$\bar E^1_{0,2n}\to H_{2n}(\g,\h, L_\mu\T L_{\la_1}^*\T L_{\la_1}^*)$$
is trivial and so $\bar E^\infty_{0,2n}=0$. Recall that
the differential $d_r$ acts from $\bar E^r_{p,q}$ to $\bar E^r_{p-r,q+r-1}$. 
This gives
$$\bar E^r_{p,q}=\bar E^1_{p,q} \text{ for } p+q\le 2n.$$
Because of 
$$H_{2n}(\g,\h, L_\mu\T L_{\la_1}^*\T L_{\la_1}^*)\simeq  \bar E^1_{n,n}$$
we have $\bar E^1_{n,n}=\bar E^\infty_{n,n}$. So 
the differential $d_n: \bar E^n_{n,n}\to \bar E^n_{0,2n-1}$ is trivial.
Therefore 
$$0=\bar E^\infty_{0,2n-1}=\bar E^n_{0,2n-1}=\bar E^1_{0,2n-1}$$
($\bar E^\infty_{0,2n-1}=0$ because by induction assumption we know that
$\bar E^\infty_{p,q}=0$ for $p+q=2n-1$, $(p,q)\ne (2n-1,0)$).
The equality $\bar E^1_{0,2n-1}=0$ gives $\bar E^1_{1,2n}=0$. Therefore,
$$0=\bar E^\infty_{0,2n}=\bar E^1_{0,2n}.$$
Theorem is proved.               
\end{proof}

\begin{cor}
For any $\la_1\in P^+_{k_1}$, $\la_2\in P^+_{k_2}$, $\mu\in P^+_{k_1+k_2}$
and $n>0$
$$H_n(\n_-, L_\mu\T L_{\la_1}^*)^{\la_2'}=0.$$ 
\end{cor}

\begin{theo}
Let $\la_1\in P^+_{k_1}$, $\la_2\in P^+_{k_2}$, $\mu\in P^+_{k_1+k_2}$. Then
\begin{equation}
\label{mform}
c_{\la_1'\la_2'}^{\mu'}(q^{-1})
=\sum_{p\ge 0} (-1)^p \sum_{\substack{w\in W\\ l(w)=p}}  q^{(w*\mu)d} 
\ch_{q^{-1}} (L_{\la_1})^{(w*\mu-\la_2)'}.
\end{equation}
Another expressions for branching functions can be obtained by interchanging
$\la_1$ and $\la_2$ in the above expression.
\end{theo}
\begin{proof}
Consider homology $H_n(\n_-, L_\mu\T L_{\la_1}^*)^{\la'_2}$. Because of 
Corollary $\ref{H0}$ and Theorem $\ref{main}$ we know that
$$\sum_{p\ge 0} (-1)^p \ch_q H_n(\n_-, L_\mu\T L_{\la_1}^*)^{\la'_2}=
c_{\la_1'\la_2'}^{\mu'}(q^{-1}).$$
Using the BGG-resolution of $L_\mu$ we obtain that the following complex
counts $H_n(\n_-, L_\mu\T L_{\la_1}^*)^{\la'_2}$:
\begin{equation}
\label{Eu}
\ldots\to D_2\to D_1\to D_0\to 0,\ D_p=\bigoplus_{w:\ l(w)=p} 
(\C\T_{\U(\n_-)} (M_{w*\mu}\T L_{\la_1}^*))^{\la'_2}.
\end{equation}
We note that 
$$(\C\T_{\U(\n_-)} (M_{w*\mu}\T L_{\la_1}^*))^{\la'_2}\simeq 
(\C_{w*\mu}\T L_{\la_1}^*)^{\la'_2}$$
and so we have  
$$\ch_q (\C\T_{\U(\n_-)} (M_{w*\mu}\T L_{\la_1}^*))^{\la'_2}=
q^{(w*\mu)d}\ch_q (L_{\la_1}^*)^{(-w*\mu + \la_2)'}.$$
Therefore the Euler characteristics of the complex $(\ref{Eu})$ is given by
$$\sum_{p\ge 0} (-1)^p \sum_{w:\ l(w)=p}  q^{(w*\mu)d} 
\ch_{q^{-1}} (L_{\la_1})^{(w*\mu-\la_2)'}.$$
But the Euler characteristics of $(\ref{Eu})$ coincides with the sum
$$\sum_{n\ge 0} (-1)^n \ch_q H_n(\n_-, L_\mu\T L_{\la_1}^*)^{\la_2'}.$$
Theorem is proved. 
\end{proof}

\section{The $\slth$ case}
In this section we specialize formulas $(\ref{genhom})$ and $(\ref{mform})$ 
to the case $\g=\slth$.

Let $h$ be the standard generator of the Cartan subalgebra of $\slt$.
Then $\h$ is spanned by $h_0=h\T 1$, $K$ and $d$. 
Define $(i,k,m)\in\h^*$ by  
$$(i,k,m)h_0=i,\    (i,k,m)K=k, \ (i,k,m)d=m.$$
Let $(i,k)=(i,k,m)|_{\h'}$.
Note that 
$$P_k^{'+}=\{(i,k)\in {\h'}^*: i,k\in \Z_{\ge 0}, i\le k\}.$$
We denote by $L_{i,k}$  ($0\le i\le k$) the highest weight irreducible 
representation of $\slth'$ with highest weight $(i,k)$. 
For fixed levels $k_1,k_2$ let $c_{i_1i_2}^j(q)$  be the corresponding 
branching functions defined as 
characters of $C_{i_1i_2}^j$:
$$L_{i_1,k_2}\T L_{i_2,k_2}=\bigoplus_{j=0}^{k_1+k_2} C_{i_1i_2}^j \T 
L_{j,k_1+k_2}.$$

Recall that for any $s>0$ there exist two elements 
$w_{s,1},w_{s,2}\in W$ with $l(w_{s,i})=s$. In addition 
\begin{gather}
\label{W}
w_{2n,1} * (i,k,m)=(i+2n(k+2),k,m+n(n(k+2)+i+1)),\\  \nonumber
w_{2n,2} * (i,k,m)=(i-2n(k+2),k,m+n(n(k+2)-i-1)),\\   \nonumber
w_{2n-1,1} * (i,k,m)=(-i-2+2n(k+2),k,m+n(n(k+2)-i-1)),\\ \nonumber
w_{2n+1,2} * (i,k,m)=(-i-2-2n(k+2),k,m+n(n(k+2)+i+1)).
\end{gather}

Let $V^l$ be the eigenspace of the operator $h\in\slt$ with
an eigenvalue $l$.  
\begin{prop}
\begin{equation}
c_{i_1i_2}^j(q)= 
\sum_{p\in\Z} q^{p^2(k_1+2)+p(i_1+1)}
\left(\ch_q L_{i_2,k_2}^{2p(k_1+2)-j+i_1}-\ch_q L_{i_2,k_2}^{2p(k_1+2)+j+i_1+2}
\right).                                       
\end{equation}
\end{prop}
\begin{proof}
Follows from $(\ref{genhom})$ and $(\ref{W})$.
\end{proof}

\begin{prop}
\begin{multline}
\label{second}
c_{i_1i_2}^j(q)=\\ 
\sum_{p\in\Z}
q^{-(k_1+k_2+2)p^2-(j+1)p}\bigl(\ch_q L_{i_1,k_1}^{2(k_1+k_2+2)p+j-i_2}-
\ch_q L_{i_1,k_1}^{2(k_1+k_2+2)p+j+i_2+2}\bigr).
\end{multline}
\end{prop}
\begin{proof}
Follows from $(\ref{mform})$ and $(\ref{W})$.
\end{proof}

\begin{rem}
Note that formula $(\ref{second})$ coincides with bosonic formula from \cite{F}.
\end{rem}

\newcounter{a}
\setcounter{a}{1}


\begin{thebibliography}{99}

\bibitem[BNY]{BNY}
J. Bagger, D. Nemeschansky, S. Yankielowicz,
{\it Virasoro algebras with central charge $c>1$},
Phys. Rev. Lett. {\bf 60} (1988), no. 5, 389-392.


\bibitem[BGG]{BGG}
I.N. Bernstein, I.M. Gel'fand, S.I. Gel'fand, {\it Differential operators on
the base affine space and a study of $\g$-modules}, in: Lie groups and
their representations (I.M. Gelfand ed.), Summer school of the Bolyai Janos
Math. Soc., Halsted Press, 1975, 21-64.

\bibitem[DJKMO]{DJKMO}
E. Date, M. Jimbo, A. Kuniba, T. Miwa, M. Okado,
{\it Exactly solvable SOS models: local heights probabilities and theta
function identities}, Nucl. Phys. B {\bf 290} (1987), no. 2, 231-273.

E. Date, M. Jimbo, A. Kuniba, T. Miwa, M. Okado,
{\it Exactly solvable SOS models \Roman{a}: proof of star-triangle relation
and combinatorial identities}, Adv. Stud. in Pure Math.
{\bf 16} (1988), 17-122.

\bibitem[DMS]{DMS}
P. Di Francesco, P. Mathieu and D. S{\'e}n{\'e}chal, {\it Conformal
field theory}, Springer GTCP, New York, 1997.


\bibitem[F]{F} 
E. Feigin, {\it Infinite fusion products and $\slth$ cosets}, preprint
2006,\\ http://xxx.lanl.gov/abs/math.QA/0603226.

\bibitem[FOW]{FOW} 
O. Foda, M. Okado, O. Warnaar, {\it A proof of polynomial identities of
type $\widehat{\mathfrak{sl}(n)_1}\T \widehat{\mathfrak{sl}(n)_1}/
\widehat{\mathfrak{sl}(n)_2}$}, J. Math. Phys. {\bf 37} (1996), 965-986.

\bibitem[FF]{FF}
B. Feigin, E. Feigin, {\it Homological realization of restricted Kostka
polynomials}, Int. Math. Res. Not. 2005, no. 33, pp. 1997-2029.


\bibitem[FFJMT]{FFJMT}
B. Feigin, E. Feigin, M. Jimbo, T. Miwa, Y. Takeyama,
{\it A $\phi_{1,3}$-filtration for Virasoro minimal series $M(p,p')$ with
$1<p'/p<2$}, preprint 2006,\\ http://xxx.lanl.gov/abs/math.QA/0603070.


\bibitem[GKO]{GKO}
P. Goddard, A. Kent, D. Olive, {\it Unitary representations of the Virasoro
and super-Virasoro algebras}, Comm. Math. Phys., {\bf 103} (1986), no.1,
105-119.

\bibitem[GL]{GL}
H. Garland, J. Lepowsky, {\it Lie algebra homology and the Macdonald-Kac
formulas}, Invent. Math. {\bf 34} (1976), 37-76.


\bibitem[Kac]{Kac}
V. Kac, {\it Infinite dimensional Lie algebras}, 3rd ed.,
Cambridge University Press, Cambridge, 1990.

\bibitem[KMQ]{KMQ}
D. Kastor, E. Martinec, Z. Qiu, {\it Current algebra and conformal decrete
series}, Phys. Lett. B{\bf 200} (1988), no. 4, 434-440.

\bibitem[Kum]{Kum} 
S. Kumar, {\it Kac-Moody groups, their flag varieties and representation 
theory}, Progressin Mathematics, vol. 204, Birkhauser Boston, Massachusetts, 
2002.


\bibitem[R]{R} F. Ravanini, {\it An infinite class of new conformal field
theories with extended algebras},  Mod. Phys. Lett. A{\bf 3} (1988), no. 4,
397-412.


\bibitem[RC]{RC}
A. Rocha-Caridi, {\it Vacuum vector representations of the Virasoro algebra},
Vertex Operators in Mathematical Physics (Berkeley, Calif., 1983),
Math. Sci. Res. Inst.  Publ., vol. 3, Springer, New York, 1985, pp. 451-473.


\bibitem[S1]{S1}
A. Schilling, {\it Multinomials and polynomial bosonic forms for the branching
functions of the
$\widehat{\mathfrak{su}}(2)_M\times \widehat{\mathfrak{su}}(2)_N/
\widehat{\mathfrak{su}}(2)_{M+N}$ conformal coset models},
Nucl. Phys. B{\bf 467} (1996), 247-271.


\bibitem[S2]{S2}
A. Schilling, {\it Polynomial fermionic forms for the branching functions
of the rational coset conformal field theories
$\widehat{\mathfrak{su}}(2)_M\times \widehat{\mathfrak{su}}(2)_N/
\widehat{\mathfrak{su}}(2)_{M+N}$}, Nucl. Phys. B{\bf 459} (1996), 393-436.


\bibitem[SS]{SS}
A. Schilling, M. Shimozono, {\it Fermionic formulas for level-restricted
generalized Kostka polynomials and coset branching functions}, Comm.
Math. Phys. {\bf 220} (2001), no. 1, 105-164.



\end{thebibliography}
\end{document}